\documentclass[12pt]{amsart}

\usepackage{fullpage}

\usepackage{amsmath}
\usepackage{amssymb}
\usepackage{amsthm}
\usepackage{graphicx}
\newtheorem{theorem}{Theorem}[section]
\newtheorem{lemma}{Lemma}[section]

\theoremstyle{definition}
\newtheorem{definition}[theorem]{Definition}

\title{Unbounded regions of Infinitely Logconcave Sequences}

\author[David Uminsky and Karen Yeats]{David Uminsky and Karen Yeats \\
Department of Mathematics and Statistics\\
Boston University, Boston, USA}

\date{March 20, 2007}
\thanks{The first author is
    partially supported by NSF grant DMS-0405724.  Thanks to Cameron Morland
    for making better figures and to the referee for a very close reading.}

\begin{document}

\begin{abstract}
We study the properties of a logconcavity operator on a symmetric,
unimodal subset of finite sequences. In doing so we are able to prove
that there is a large unbounded region in this subset that is
$\infty$-logconcave.  This problem was motivated by the conjecture of
Boros and Moll in \cite{bomo} that the binomial coefficients are
$\infty$-logconcave. 
\end{abstract}

\maketitle

\section{Introduction}
In this paper we study the asymptotic behavior of the logconcavity operator on finite sequences. Before we can state the problem we will need a few definitions.

 We say that a sequence $\{c_0,c_1,...c_n\}$ is \textbf{$1$-logconcave}
 (or logconcave) if $c_i\ge 
0$ for
$0\le i \le n$ and $c_i^2 -c_{i+1}c_{i-1} \ge 0$ for
 $1\le i \le n-1$. 

We can extend this idea of logconcave as follows: 
Since $\{c_i\}$ is a finite sequence of length $n$ we define $c_i =0$ for $i<0$ and $i>n+1$, then define the operator
\begin{equation}
\mathcal{L}\{c_i\} = \{c_i^2-c_{i-1}c_{i+1} \}. 
\end{equation}  

If $\{c_i\}$ is logconcave then $\mathcal{L}\{c_i\}$ is a new non-negative
sequence. We now define a sequence $\{c_i\}$ to be
\textbf{$\infty$-logconcave} if $\mathcal{L}^k\{c_i\}$ is a non-negative
sequence for all $k\ge 1$. 

While studying a new class of integrals related to Ramanujan's Master
Theorem, Boros and Moll proposed that a particular family of finite
sequences of coefficients $\{d_l(m)\}$  were $\infty$-logconcave. Boros
and Moll then point out that showing that the binomial coefficients
are $\infty$-logconcave (project 7.9.3 in \cite{bomo}) would go a long
way in showing the sequence $\{d_l(m)\}$ is $\infty$-logconcave.
Kauers and Paule in \cite{kp} show that the $\{d_l(m)\}$ are $1$-logconcave.
These conjectures motivated us to investigate the operator $\mathcal{L}$ on the space of finite sequences.  

%
%
Numerical experiments suggest that the binomial coefficients are
$\infty$-logconcave.  Moreover, many sequences ``near'' the binomial
sequence also appear to be $\infty$-logconcave.  These numerics led us to take an alternative approach.  We begin to study the properties of $\mathcal{L}$ on the subset of finite sequences of the forms
\[
  \{\dots,0,0,1, x_0, x_1, \dots, x_n, \dots, x_1, x_0, 1, 0, 0,
  \dots\}
\]
\[
  \{\dots,0,0,1, x_0, x_1, \dots, x_n, x_n, \dots, x_1, x_0, 1, 0, 0,
  \dots\}
\]

 We will refer to the first sequence above as the odd case and to the
 second sequence as the even case because of the repetition of the middle term in the sequence.  Notice that all the binomial coefficients belong to one of the above cases.

Our approach to the problem is the following: for a given sequence of the form above
of length $2n+3$ or $2n+4$ we  analyze the dynamics of $\mathcal{L}$
on the subset of Euclidean space $\mathbb{R}^{n}$ with all non-negative
coordinates.  This differs from the approach
of Moll in \cite{m}.

Our main result in this paper is to show that there is a large
unbounded region $\mathcal{R}$ in this orthant that contains only
$\infty$-logconcave sequences.  Moreover $\mathcal{R}$ can act like a
trapping region for $\infty$-logconcave sequences, i.e., sequences not
starting in $\mathcal{R}$ can land in $\mathcal{R}$ after a number of
iterates of $\mathcal{L}$. 


The paper is organized as follows: In section \ref{LowDim} we present some of the simple cases along with some numerical evidence. The general arguments are presented in detail for the even case in section \ref{evencase} and the odd case is briefly covered in similar fashion in section \ref{oddcase}.

\section{Low Dimensional Cases}\label{LowDim}

\subsection{The one dimensional cases, $\{1,x,1\}$ and $\{1,x,x,1\}$}

For these two cases the underlying dynamics of $\mathcal{L}$ is rather easy to compute explicitly. We first consider the sequence $\{1,x,1\}$.

\begin{theorem}\label{1Dodd}
$\mathcal{L} \{1,x,1\} = \{1,x^2-1,1\}$. Thus the positive fixed point for this
sequence is $x = \frac{1+\sqrt{5}}{2}$.  Moreover, if $x\ge \frac{1+\sqrt{5}}{2}$ then our sequence is $\infty$-logconcave and not otherwise.
\end{theorem}
\begin{proof}
A simple calculation shows that $x = \frac{1+\sqrt{5}}{2}$ is a fixed
point of $\mathcal{L}$ for the sequence  $\{1,x,1\}$.  Moreover it is
easy to see that if if $x > \frac{1+\sqrt{5}}{2}$ then $x$ grows under
the iterates of $\mathcal{L}$ and hence is always positive.  It is also easy to see that the interval $[1,\frac{1+\sqrt{5}}{2}]$ is mapped over the interval $[0,\frac{1+\sqrt{5}}{2}]$ monotonically so that any values of $x \in [1, \frac{1+\sqrt{5}}{2})$ are eventually mapped below $x=1$. Therefore our sequence is no longer unimodal and, hence, not logconcave.
\end{proof}

Notice that the binomial sequence $\{1,2,1\}$ lies securely in this region thus we have shown that $\{1,2,1\}$ is $\infty$-logconcave. The sequence $\{1,x,x,1\}$ is handled in a similar fashion.

\begin{theorem}\label{1Deven}
$\mathcal{L} \{1,x,x,1\} = \{1,x^2-x,x^2-x,1\}$, thus the fixed point for this sequence is $x = 2$. If $x\ge 2$ then our sequence is $\infty$-logconcave and not otherwise.

\end{theorem}
\begin{proof}
Nearly identical to the one above.
\end{proof}

Notice that the key to the easy theorems above is finding the fixed
points of $\mathcal{L}$ in our underlying Euclidean space and then
monotonicity leads us to the rest.  Fixed ``points'' will no longer remain the key in the general argument but as we will see in the 2-D cases below we will have hypersurfaces that bound open regions of $\infty$-logconcavity.

\subsection{The 2-D cases, $\{1,x,y,x,1\}$ and $\{1,x,y,y,x,1\}$}

The 2-D cases are more complicated than the 1-D cases but they give better insight into how we might hope to find ``regions" of $\infty$-logconcavity. To better see the general techniques to finding such regions, we first focus on the even case.


If one is to investigate this question numerically one can compute the following picture.

\begin{figure}[htbp]
\begin{center}
\includegraphics[width=0.35\textwidth ]{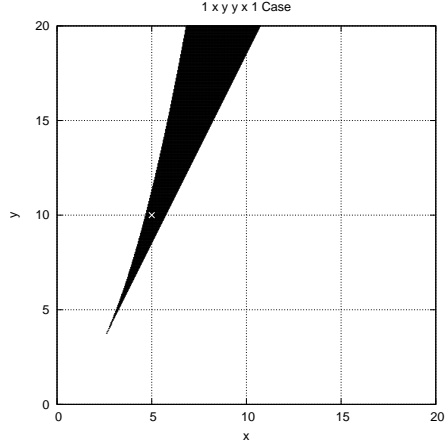}
\caption{The filled region is the numerical region of $\infty$-logconcavity for the 2-D even case. The X indicates the position of the binomial coefficient.}
\end{center}
\end{figure}

The first thing to notice is that the binomial coefficient $x=5,y=10$
is in the numerical region of $\infty$-logconcavity.  This picture also suggests that there is an
$\infty$-logconcave region  bounded
away from the origin, below by some line and above by some curve.  This picture is remarkably
stable.  The boundary points in the 1-D cases have now been replaced by
curves.

For the case $\{1,x,y,x,1\}$ the results are similar. The numerical picture looks as follows.

\begin{figure}[htbp]
\begin{center}
\includegraphics[width=0.35\textwidth ]{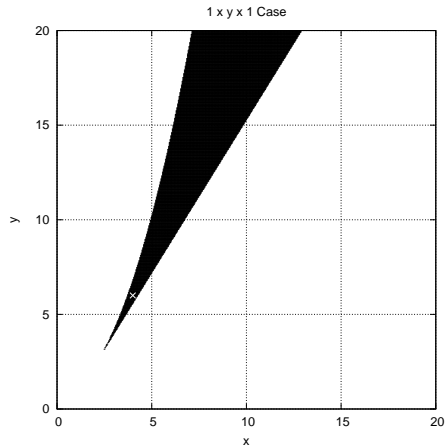}
\caption{The filled region is the numerical region of $\infty$-logconcavity for the 2-D odd case.  The X indicates the position of the binomial coefficient. }
\end{center}
\end{figure}

Again, we notice that the binomial coefficient $x=4,y=6$ is in the
numerical region of $\infty$-logconcavity, which is quite
encouraging. This picture also suggests that there is a region of
$\infty$-logconcavity bounded away from the origin.  It is important to point out that the regions in both cases are different with the odd case containing a wider region.

\subsection{Note on a 3-D case}\label{3d}

If one is to investigate the even 3-D case $\{1, x, y, z,
z, y, x, 1\}$, we would arrive at the following boundary
hypersurfaces.  
\begin{figure}[htbp]
\begin{center}
\input{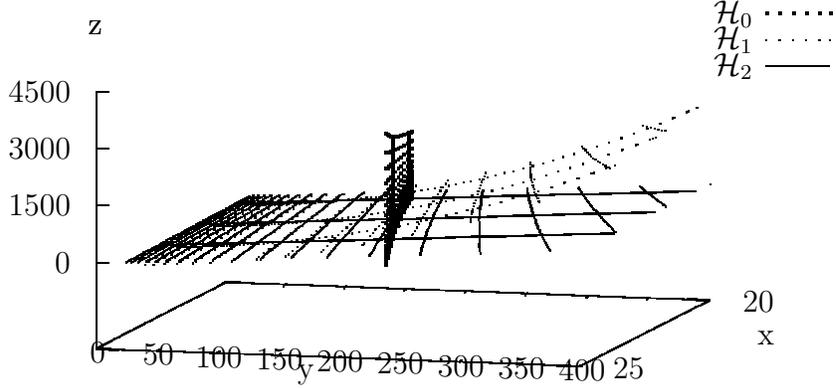}
\caption{The region of $\infty$-logconcavity for the 3-D even case. From the picture above the region of interest is bounded on the left by the ``vertical" plane ($\mathcal{H}_0$), below by the ``horizontal" plane ($\mathcal{H}_2$), and above by the curved surface ($\mathcal{H}_1$).}\label{3DPlot}
\end{center}
\end{figure}

In this case the binomial coefficient sequence $x=7, y=21, z=35$ is
not in the region.  However its first iterate $x=28, y=196, z=490$ is in the region of interest.

It is the observation of these hypersurfaces as boundaries of the
$\infty$-logconcave region that is important.  As it turns out, we can
construct these boundaries for arbitrarily long finite sequences as
will be shown in section \ref{evencase} for the even cases and
section \ref{oddcase} for the odd cases. While the 1-D and 2-D cases
seem to have an ``exclusive'' region of $\infty$-logconcavity, this is
not true in general as we saw in the 3-D case. 

\section{The General Case for Even Length}\label{evencase}
In this section we will prove the existence of an unbounded region of
infinite logconcavity for even length symmetric sequences.  Before
beginning in earnest, let us mention the main steps.  By looking at the leading order behavior we find hypersurfaces
bounding the region of infinite logconcavity.  We proceed to show the
region is nonempty and unbounded by an explicit example, which we can
control by matching all but one of the coordinates of the example with
each hypersurface in turn.  Then we show that sequences within the
region are indeed $\infty$-logconcave by again matching with the
hypersurfaces in turn and iterating.  The technical key to controlling
the iteration is understanding the effect of increasing one coordinate
while decreasing all the others.

\subsection{Leading order behavior}
Consider the sequence of length $2n+4$:
\[
 s = \{1, a_0x, a_1x^{1+d_1}, a_2x^{1+d_1+d_2}, \ldots,
 a_nx^{1+d_1+\cdots+d_n}, a_nx^{1+d_1+\cdots+d_n}, \ldots, a_0x, 1\}.
\]
For the moment we are interested in the leading terms of elements of
$\mathcal{L}(s)$ viewed as polynomials
in $x$.  We will restrict ourselves to values of the $d_i$ for which
$a_i^2x^{2(1+d_1+\cdots+d_i)}$ contributes to the leading term of the
corresponding element in the first iteration
\begin{equation}\label{iterate}
\begin{split}
 \{1, x(a_0^2x - a_1x^{d_1}), x^{2+d_1}(a_1^2x^{d_1} - a_2a_0x^{d_2}),
 x^{2+2d_1+d_2}(a_2^2x^{d_2} - a_1a_3x^{d_3}), \ldots, \\
 x^{2+2d_1+\cdots+2d_{n-1}+d_n}(a_n^2x^{d_n} -
 a_na_{n-1}), x^{2+2d_1+\cdots+2d_{n-1}+d_n}(a_n^2x^{d_n} -
 a_na_{n-1}), \ldots, 1 \}
\end{split}
\end{equation}
and restrict ourselves to values of $a_i > 0$ which give that the leading terms of
$\mathcal{L}(s)$ have the same form as $s$ itself for some new $x$.

Using (\ref{iterate}), the leading term condition is equivalent to
\[
  0 \leq d_n \leq d_{n-1} \leq \cdots \leq d_1 \leq 1
\]
which we can view as defining a simplex.  The values of $a_i$ may then
be determined for each face by solving the systems of equations
arising from matching the coefficients of the leading terms in
$\mathcal{L}(s)$ with the 
coefficients of $s$.  

Of greatest interest are the $(n-1)$-faces of the simplex since they define the
boundaries of what will be our open region of convergence.  The $(n-1)$-faces are defined by $d_1=1$, $d_j=d_{j+1}$ for $0 < j < n$, and
$d_n=0$; in all cases the unspecified $d_i$ are distinct and strictly
between $0$ and $1$.

For $d_1=1$ the leading terms of $\mathcal{L}(s)$ are
\[
  \{1, x^2(a_0^2 - a_1), x^{4}a_1^2, x^{4+2d_2}a_2^2, \ldots,
  x^{4+2d_2+\cdots+2d_n}a_{n}^2, x^{4+2d_2+\cdots+2d_n}a_{n}^2, \ldots, 1\}
\]
so we are led to the system
\begin{align*}
  a_0^2 - a_1 & = a_0 \\
  a_1^2 & = a_1 \\
  & \vdots \\
  a_n^2 & = a_n
\end{align*}
We are interested in positive solutions so $a_i = 1$ for $0 < i \leq
n$ and $a_0 = (1+\sqrt{5})/2$.

For $d_j = d_{j+1}$ the leading terms of $\mathcal{L}(s)$ are
\[
\begin{split}
  \{1, x^2 a_0^2, x^{2+2d_1}a_1^2, \ldots,
  x^{2+2d_1+\cdots+2d_j}(a_j^2 - a_{j-1}a_{j+1}),
  x^{2+2d_1+\cdots+4d_j}a_{j+1}^2, \ldots, \\
  x^{2+2d_1+\cdots+4d_j+\cdots+2d_n}a_{n}^2, x^{2+2d_1+\cdots+4d_j+\cdots+2d_n}a_{n}^2, \cdots, 1\}
\end{split}
\]
so we are led to the system
\begin{align*}
  a_0^2 & = a_0 \\
  & \vdots \\
  a_{j-1}^2 & = a_{j-1} \\
  a_j^2 -a_{j-1}a_{j+1} & = a_j \\
  a_{j+1}^2 & = a_{j+1} \\
  & \vdots \\
  a_n^2 & = a_n
\end{align*}
which has unique positive solution $a_i = 1$ for $i \neq j$ and $a_j =
(1+\sqrt{5})/2$. 

Finally for $d_n=0$ the leading terms of $\mathcal{L}(s)$ are
\[
\begin{split}
  \{1, x^2a_0^2, x^{2+2d_1}a_1^2, \ldots,
  x^{2+2d_1+\cdots+2d_{n-1}}a_{n-1}^2, \\
  x^{2+2d_1+\cdots+2d_{n-1}}(a_{n}^2-a_na_{n-1}),
  x^{2+2d_1+\cdots+2d_{n-1}}(a_{n}^2-a_na_{n-1}), \ldots, 1\}
\end{split}
\]
so we are led to the system
\begin{align*}
  a_0^2 & = a_0 \\
  & \vdots \\
  a_{n-1}^2 & = a_{n-1} \\
  a_n^2 -a_na_{n-1} & = a_n
\end{align*}
which has unique positive solution $a_i = 1$ for $0 \leq i <
n$ and $a_n = 2$.

\subsection{Interior}
In the region of $\mathbb{R}^{n+1}$ where the coordinates are all
 positive and increasing,  consider the following parametrically defined
 hypersurfaces:
\begin{align*}
  \mathcal{H}_0 & = \left\{\left(\frac{1+\sqrt{5}}{2}x, x^2,
  x^{2+d_2}, \ldots, x^{2+d_2+\cdots+d_n}\right) : 1 \leq x, 1 > d_2 >
  \cdots > d_n > 0\right\} \\
  \mathcal{H}_j & = \Bigg\{\left(x, x^{1+d_1}, \ldots,
  \frac{1+\sqrt{5}}{2}x^{1+d_1+\cdots+d_j}, x^{1+d_1+\cdots+2d_j},
  \ldots, x^{1+d_1+\cdots+2d_j+\cdots+d_n}\right) \\
  & \qquad : 1 \leq x, 1 > d_1
  > \cdots > d_j > d_{j+2} > \cdots > d_n > 0 \Bigg\} \\
  \mathcal{H}_n & = \Bigg\{\left(x, x^{1+d_1}, \ldots,
      x^{1+d_1+\cdots+d_{n-1}}, 2x^{1+d_1+\cdots+d_{n-1}}\right)
     : 1 \leq
  x, 1 > d_1 > \cdots > d_{n-1} > 0\Bigg\}
\end{align*}
for $0 < j < n$.  These are precisely the results of the leading order
analysis of the previous subsection.

Let $\mathcal{R}$ be the region with positive increasing coordinates
defined as greater in the $i$th coordinate than $\mathcal{H}_i$.  For
example in the 3-D case handled in section \ref{3d}, figure \ref{3DPlot}, the region in question is above
$\mathcal{H}_2$, below $\mathcal{H}_1$ and to the right of
$\mathcal{H}_0$.

We say a sequence $\{1, x_0, \ldots, x_n, x_n, \ldots, x_0, 1\}$ is in
$\mathcal{R}$ if $(x_0, \ldots, x_n) \in \mathcal{R}$.

Before we discuss $\mathcal{R}$ further we must first recall that the $n^{th}$ 
{\bf triangular number}, $\tilde{T}(n)$,  is defined as $$\tilde{T}(n) = \tilde{T}(n-1) +n,$$ with $\tilde{T}(0)=0$.  The first few elements of the sequence are $0,1, 3, 6, 10, 15, 21, 28, 36, 45, \ldots $. We will need the following lemma about triangular numbers.
\begin{lemma}\label{triangular}
Define $T(n) \equiv 2\tilde{T}(n)$ for $n \geq 0$. Then $T(n)$ satisfies the following:
\begin{enumerate}
\item $T(0) -\frac{T(1)}{2} =-1$
\item $T(n+1) = 2T(n) -T(n-1) +2$ 
\end{enumerate}
\end{lemma}

\begin{proof}
(1) is trivial. (2) follows since $T(n) = n(n+1)$.
\end{proof}

We also need another straightforward result.
\begin{lemma}\label{1log}
Suppose $x>0$ and 
\[
s = \{1, x, x^{1+d_1}, x^{1+d_1+d_2}, \cdots,
x^{1+d_1+\cdots+d_n}, x^{1+d_1+\cdots+d_n}, \cdots , x, 1\}.
\]
  Then
$s$ is $1$-logconcave iff $1 \geq d_1 \geq \cdots \geq d_n \geq 0$
with strict inequalities in the logconcavity condition iff $1 > d_1 > \cdots > d_n > 0$.
\end{lemma}

\begin{proof}
  Compute
  $x^2 \geq x^{1+d_1} \Leftrightarrow 1 \geq d_1$ and $x^2 > x^{1+d_1}
  \Leftrightarrow 1 > d_1$.  If $0 < j <
  n$ then  $x^{2+2d_1+\cdots+2d_{j}} \geq x^{2+2d_1 + \cdots + 2d_{j-1}
  + d_{j}+ d_{j+1}} \Leftrightarrow x^{d_j} \geq x^{d_{j+1}}
  \Leftrightarrow d_j \geq d_{j+1}$ and likewise with strict inequalities.  Finally
  $x^{2+2d_1+\cdots+2d_n} \geq x^{2+2d_1+\cdots+2d_{n-1} + d_n}
  \Leftrightarrow x^{d_n} \geq 1 \Leftrightarrow d_n \geq 0$ and likewise
  with strict inequalities.
\end{proof}

We are now ready to prove some important properties of $\mathcal{R}$.

\begin{lemma}
  $\mathcal{R}$ is nonempty and unbounded.
\end{lemma}

\begin{proof}
Let $\{1, x_0, \ldots,
x_n, x_n, \ldots, x_0, 1\}$ be any $1$-logconcave sequence with $x_0 >
0$, for
instance the binomial sequence of appropriate length. Also, choose $C$ such that $0<C < \frac{2}{1+\sqrt{5}}$ and
consider the following sequence:

$$
s=\Bigg\{1, C^{T(0)}ax_0,C^{T(1)} a^2 x_1,
 C^{T(2)} a^3 x_2, \ldots, 
 C^{T(n)} a^{n+1}x_n, 
  C^{T(n)} a^{n+1}x_n, \ldots,
  1 \Bigg\}
$$

for $a > 2 C^{T(n-1)-T(n)}$.




Notice that $a$ is dependent on $n$ which is not a problem since $n$ is fixed.

It is clear that $s$ is $1$-logconcave and, moreover, the inequalities are strict since
\[
  C^{2T(0)}a^2x_0^2 = a^2x_0^2 \geq
  a^2x_1 > C^{T(1)}a^2x_1,
\]
for $0 < j < n$
\begin{align*}
  C^{2T(j)}a^{2j+2}x_j^2 &\geq
  C^{2T(j)}a^{2j+2}x_{j-1}
  x_{j+1} \\ & > C^{T(j-1)}a^{j}x_{j-1}
  C^{T(j+1)}a^{j+2}x_{j+1}
\end{align*}
by (2) of Lemma \ref{triangular} and the fact that $C < 2/(1+\sqrt{5}) <
1$, and 
\begin{align*}
  C^{T(n)}a^{n+1}x_n &\geq
  C^{T(n)}a^{n+1}x_{n-1} \\ 
  & > C^{T(n-1)}a^{n}x_{n-1}.
\end{align*}

Define $\tilde{x} = ax_0 > 0$, define $\tilde{d}_1$ so that
$\tilde{x}^{1+\tilde{d}_1} = C^{T(1)}a^2x_1$, and
continue recursively so that, for $0 < j \leq n$, $\tilde{d}_j$ is
defined so that $\tilde{x}^{1+\tilde{d}_1 +\cdots + \tilde{d}_j} =
C^{T(j)}a^{j+1}x_j$.  By Lemma \ref{1log}, $1 > \tilde{d}_1 >
\cdots > \tilde{d}_n > 0$.

Let us next consider each $\mathcal{H}_j$ in turn. For $0 < j < n$ choose $x=\tilde{x}$, $d_i = \tilde{d}_i$ for
$i<j$, $d_j = (\tilde{d}_j + \tilde{d}_{j+1})/2$, and $d_i =
\tilde{d}_i$ for $i > j+1$.  Consequently $1 > d_1 > \cdots > d_j >
d_{j+2} > \cdots > d_n > 0$ and 
these choices match all the coordinates of $s$ with the corresponding coordinates of
$\mathcal{H}_i$ except possibly for the
$j$th.  But $x^{1+d_1+\cdots+2d_{j}}/x^{1+d_1+\cdots+d_{j-1}} =
x^{2d_{j}}$, so 
\begin{align*}
 Cx^{1+d_1+\cdots+d_j} & =
 C x^{1+d_1+\cdots+d_{j-1}}\sqrt{\frac{x^{1+d_1+\cdots+2d_{j}}}{x^{1+d_1+\cdots+d_{j-1}}}} \\
 & =
 C\sqrt{x^{1+d_1+\cdots+d_{j-1}}x^{1+d_1+\cdots+2d_{j}}}.
\end{align*}
Comparing with $s$ we have that
\begin{align*}
  C^{T(j)} a^{j+1}x_j
  & \geq C^{T(j)}
  a^{j+1}\sqrt{x_{j-1}x_{j+1}} \\
  & =
  \sqrt{C^{2T(j)-T(j+1)-T(j-1)}
  C^{T(j-1)}a^{j}x_{j-1}   
  C^{T(j+1)}a^{j+2}x_{j+1}} \\
  & = C^{-1}\sqrt{x^{1+d_1+\cdots+d_{j-1}}x^{1+d_1+\cdots+2d_{j}}} \\
  & > \frac{1+\sqrt{5}}{2}x^{1+d_1+\cdots+d_j}
\end{align*}
where the fourth line follows from $(2)$ of Lemma \ref{triangular}, thus $s$ is on the correct side of $\mathcal{H}_j$.

For $\mathcal{H}_0$ choose $x = \sqrt{C^{T(1)}a^2x_1} =
\tilde{x}^{(1+\tilde{d}_1)/2} > 0$ and, for $2 \leq j \leq n$, $d_j
= 2\tilde{d_j}/(1+\tilde{d_1})$.  Consequently $1 > d_2 > \cdots > d_n
> 0$ and 
\[
C^{T(j)}a^{j+1}x_j =
\tilde{x}^{1+\tilde{d}_1+\cdots+\tilde{d}_j} = x^{2+d_2+\cdots+d_j}
\]
hence
matching all the coordinates of $s$ other than the
$0$th with $\mathcal{H}_0$.  So we check,
\[
  C^{T(0)}ax_0 \geq C^{T(0)}\sqrt{a^2x_1} = C^{T(0)}C^{-\frac{T(1)}{2}}x = C^{-1}x >\frac{1+\sqrt{5}}{2}x
\]
Thus $s$ is on the correct side of $\mathcal{H}_0$.

For $\mathcal{H}_n$ simply choose $x=\tilde{x}$ and $d_i = \tilde{d}_i$ for
$i<n$, which gives $1 > d_1 > \cdots > d_{n-1} > 0$ and matches
all the coordinates of $s$ other than the $n$th.  Then $x_n^2 \geq
x_{n-1}x_n$ giving
\begin{align*}
  C^{T(n)}a^{n+1}x_{n} 
  & \geq
  C^{T(n)}a^{n+1}x_{n-1}
  \\
  & = a C^{T(n)-T(n-1)} x^{1+d_1+\cdots +d_{n-1}} \\
  & > 2x^{1+d_1+\cdots +d_{n-1}}
\end{align*}
So $s$ is also on the correct side of $\mathcal{H}_n$.  Consequently
$s$ is in $\mathcal{R}$.  So we see that $\mathcal{R}$ is nonempty,
and, by the freedom to increase $a$, is unbounded.
\end{proof}

\begin{definition}
  Let $\mathcal{H}$ be a hypersurface in $\mathbb{R}^{n+1}$.
  We say we view $\mathcal{H}$ as a function $f : \mathbb{R}^n \rightarrow \mathbb{R}$ with
  the $j$th variable as the dependent variable if for $(x_0,
  \ldots, x_n)$ a point on $\mathcal{H}$ we have $x_j = f(x_0, \ldots,
  x_{j-1}, x_{j+1}, \ldots, x_n)$.
\end{definition}

\begin{definition}
  Let $\mathcal{H}$ be a hypersurface in $\mathbb{R}^{n+1}$.  Call it
  \textbf{$j$-monotone} if when $\mathcal{H}$ is
  viewed as a function $f : \mathbb{R}^n \rightarrow \mathbb{R}$ with
  the $j$th variable as the dependent variable then $f(y_1,
  \ldots, y_n) \geq f(z_1, \ldots, z_n)$ if $y_i \geq z_i$ for all $i$. 
\end{definition}

\begin{lemma}\label{sides}
  Let $\mathcal{H}$ be a $j$-monotone hypersurface in
  $\mathbb{R}^{n+1}$.
  Let
  $(x_0, x_1, \ldots, x_n)$ be a point on $\mathcal{H}$.  Then for
  $\epsilon_i > 0$, $\eta > 0$, 
  \[
  (x_0-\epsilon_0, 
  \ldots, x_{j-1}-\epsilon_{j-1}, x_j+\epsilon_j,
  x_{j+1}-\epsilon_{j+1}, \ldots, x_n-\epsilon_n)
  \]
  and 
  \[
  (x_0, \ldots, x_{j-1}, x_j+\eta, x_{j-1}, \ldots, x_n)
  \]
  lie on the same side of $\mathcal{H}$.
\end{lemma}

\begin{proof}
  View $\mathcal{H}$ as $f : \mathbb{R}^n \rightarrow \mathbb{R}$ with
  $x_j$ as the dependent variable.  Then 
  \begin{align*}
    & f(x_0-\epsilon_0, 
    \ldots, x_{j-1}-\epsilon_{j-1},
    x_{j+1}-\epsilon_{j+1}, \ldots, x_n-\epsilon_n) \\
    & \leq f(x_0,
    \ldots, x_{j-1},  x_{j-1}, \ldots, x_n) = x_j < x_j + \epsilon_j
  \end{align*}
  So both points lie on the side of $\mathcal{H}$ which is greater in
  the $j$th coordinate.
\end{proof}

\begin{lemma}\label{Hmonotone}
  Each of the $\mathcal{H}_j$ is $j$-monotone.
\end{lemma}

\begin{proof}
  For $\mathcal{H}_0$, $x_0$ is determined by $x_1$ and increases when
  $x_1$ increases, so $\mathcal{H}_0$ is $0$-monotone.  For
  $\mathcal{H}_n$, $x_n$ is determined by $x_{n-1}$ and increases when
  $x_{n-1}$ increases, so $\mathcal{H}_n$ is $n$-monotone.  For
  $\mathcal{H}_j$, $0< j < n$, $x_j$ is
  $(1+\sqrt{5})\sqrt{x_{j-1}x_{j+1}}/2$ which increases when either
  $x_{j-1}$ or $x_{j+1}$ increase so $\mathcal{H}_j$ is $j$-monotone.
\end{proof}

\begin{theorem}
  Any sequence in $\mathcal{R}$ is $\infty$-logconcave.
\end{theorem}

\begin{proof}
  Suppose $s = \{1, y_0, \ldots, y_n, y_n, \ldots, y_0, 1\}$ is in
  $\mathcal{R}$.  Then for any $0 < j < n$, by the definition of
  $\mathcal{R}$, we can choose $x$, 
  $\epsilon>0$, and the 
  $d_i$, $i\neq j+1$, such that 
  \[
  \begin{split}
    s = \Bigg\{1, x, \ldots, x^{1+d_1+\cdots+d_{j-1}},
      \frac{1+\sqrt{5}}{2}x^{1+d_1+\cdots+d_j}+\epsilon,
      x^{1+d_1+\cdots+2d_j}, \ldots, \\x^{1+d_1+\cdots+2d_j+\cdots+d_n},
       x^{1+d_1+\cdots+2d_j+\cdots+d_n}, \ldots, 1\Bigg\}.
  \end{split}
  \]
  Iterate to get



  \begin{eqnarray*}
    \mathcal{L}(s)& = & \Bigg\{1, x^2 - x^{1+d_1}, \ldots,  x^{2+2d_1+\cdots+2d_{j-1}}  \\
    & & -
    \frac{1+\sqrt{5}}{2}x^{2+2d_1+\cdots+2d_{j-2}+d_{j-1}+d_j} -
    \epsilon x^{1+d_1+\cdots+d_{j-2}}, \\& &
    \left(\left(\frac{1+\sqrt{5}}{2}\right)^2 -
    1\right)x^{2+2d_1+\cdots+2d_j} + 
    (1+\sqrt{5})x^{1+d_1+\cdots+d_j}\epsilon + \epsilon^2, \\ & &
    x^{2+2d_1+\cdots+4d_{j}} -
    \frac{1+\sqrt{5}}{2}x^{2+2d_1+\cdots+3d_{j}+d_{j+1}} -
    \epsilon x^{1+d_1+\cdots+2d_j+d_{j+1}}, \ldots, \\ & &
    x^{2+2d_1+\cdots+4d_j+\cdots+2d_n} -
    x^{2+2d_1+\cdots+4d_j+\cdots+d_n}, \\& &
    x^{2+2d_1+\cdots+4d_j+\cdots+2d_n} -
    x^{2+2d_1+\cdots+4d_j+\cdots+d_n}, \ldots, 1\Bigg\}
  \end{eqnarray*}

  Since $(1+\sqrt{5})/2)^2 -1 = (1+\sqrt{5})/2$ by using $x^2$ in place
  of $x$ in the definition of $\mathcal{H}_j$ and applying Lemma \ref{sides}, which is valid in view of
  Lemma \ref{Hmonotone}, we can conclude that $\mathcal{L}(s)$ is on
  the side of $\mathcal{H}_j$ which is larger in the $j$th
  coordinate. This is the same side which $s$ is on.

  Similarly for $\mathcal{H}_0$ we can choose $x$, $\epsilon>0$, and
  the $d_i$, such that
  \[
    s = \Bigg\{1, \frac{1+\sqrt{5}}{2} x + \epsilon, x^{2}, \ldots,
    x^{2+d_2+\cdots+d_n}, 
    x^{2+d_2+\cdots+d_n}, \ldots, 1\Bigg\}
   \]
  Iterate
  to get
  \[
  \begin{split}
    \mathcal{L}(s) = \Bigg\{1, \left(\left(\frac{1+\sqrt{5}}{2}\right)^2 -
    1\right) x^2 + \epsilon (1+\sqrt{5})x + \epsilon^2, x^4 -
    \frac{1+\sqrt{5}}{2}x^{3+d_2}-\epsilon x^{2+d_2}, \ldots, \\
    x^{4+2d_2+\cdots+2d_n} -
    x^{4+2d_2+\cdots+d_n}, 
    x^{4+2d_2+\cdots+2d_n} -
    x^{4+2d_2+\cdots+d_n}, \ldots, 1\Bigg\}
  \end{split}
  \]
  which, by Lemmas \ref{sides} and \ref{Hmonotone} shows
  that $\mathcal{L}(s)$ is on the same side of $\mathcal{H}_0$ as $s$
  is, as above.

  Finally for $\mathcal{H}_n$ choose $x$, $\epsilon>0$, and
  the $d_i$, such that
  \[
    s = \Bigg\{1, x, \ldots, x^{1+d_1+\cdots+d_{n-1}},
    2x^{1+d_1+\cdots+d_{n-1}} + \epsilon, 
    2x^{1+d_1+\cdots+d_{n-1}} + \epsilon, \ldots, 1\Bigg\}
   \]
  Iterate
  to get
  \[
  \begin{split}
    \mathcal{L}(s) = \Bigg\{1, x^2-x^{1+d_1}, \ldots,
    x^{2+2d_1+\cdots+2d_{n-1}} - 2x^{2+2d_1+\cdots+d_{n-1}}, \\
    (4-2)x^{2+2d_1+\cdots+2d_{n-1}} + 4\epsilon
    x^{1+d_1+\cdots+d_{n-1}} + \epsilon^2, \\
    (4-2)x^{2+2d_1+\cdots+2d_{n-1}} + 4\epsilon
    x^{1+d_1+\cdots+d_{n-1}} + \epsilon^2, 
    \ldots, 1\Bigg\}
  \end{split}
  \]
  which again by Lemmas \ref{sides} and \ref{Hmonotone} shows
  that $\mathcal{L}(s)$ is on the same side of $\mathcal{H}_n$ as $s$
  is.

  Consequently $\mathcal{L}(s)$ is in $\mathcal{R}$.  Since
  $\mathcal{R}$ is a subregion of the region of $\mathbb{R}^{n+1}$
  with positive coordinates, this implies that any sequence in
  $\mathcal{R}$ is $\infty$-logconcave.
\end{proof}

\section{The General Case for Odd Length}\label{oddcase}
\subsection{Leading order behavior}
Consider the sequence of length $2n+3$
\[
 s = \{1, a_0x, a_1x^{1+d_1}, a_2x^{1+d_1+d_2}, \ldots,
 a_nx^{1+d_1+\cdots+d_n}, \ldots, a_0x, 1\}
\]
Again we are interested in the leading terms of elements of
$\mathcal{L}(s)$ viewed as polynomials
in $x$.  We will restrict ourselves to values of the $d_i$ for which
$a_i^2x^{2(1+d_1+\cdots+d_i)}$ contributes to the leading term of the
corresponding element in the first iteration:
\begin{equation}\label{odditerate}
\begin{split}
 \{1, x(a_0^2x - a_1x^{d_1}), x^{2+d_1}(a_1^2x^{d_1} - a_2a_0x^{d_2}),
 x^{2+2d_1+d_2}(a_2^2x^{d_2} - a_1a_3x^{d_3}), \ldots, \\
 x^{2+2d_1+\cdots+2d_{n-1}}(a_n^2x^{2d_n} -
 a_{n-1}^2), \ldots, 1 \}
\end{split}
\end{equation}
and to values of $a_i > 0$ which give that the leading terms of
$\mathcal{L}(s)$ have the same form as $s$ itself for some new $x$.

Using (\ref{odditerate}), the leading term condition is equivalent to
\[
  0 \leq d_n \leq d_{n-1} \leq \cdots \leq d_1 \leq 1
\]
which we can again view as defining a simplex. 
The $(n-1)$-faces are defined by $d_1=1$, $d_j=d_{j+1}$ for $0 \leq j < n$, and
$d_n=0$; in all cases with the unspecified $d_i$ distinct and strictly
between $0$ and $1$.

For $d_n=0$ the leading terms of $\mathcal{L}(s)$ are
\[
\begin{split}
  \{1, x^2a_0^2, x^{2+2d_1}a_1^2, \ldots,
  x^{2+2d_1+\cdots+2d_{n-1}}a_{n-1}^2,
  x^{2+2d_1+\cdots+2d_{n-1}}(a_{n}^2-a_{n-1}^2), \ldots, 1\}
\end{split}
\]
so we are led to the system
\begin{align*}
  a_0^2 & = a_0 \\
  & \vdots \\
  a_{n-1}^2 & = a_{n-1} \\
  a_n^2 -a_{n-1}^2 & = a_n
\end{align*}
which has unique positive solution $a_i = 1$ for $0 \leq i <
n$ and $a_n = (1+\sqrt{5})/2$.

For $d_1=1$ and $d_j = d_{j+1}$ the systems are identical to the even case. 

\subsection{Interior}
The hypersurfaces $\mathcal{H}_j$, $0 \leq j < n$ are the same. 
\begin{align*}
  \mathcal{H}_n & = \Bigg\{\left(x, x^{1+d_1}, \ldots,
      x^{1+d_1+\cdots+d_{n-1}}, \frac{1+\sqrt{5}}{2}x^{1+d_1+\cdots+d_{n-1}}\right) \\
    & : 1 \leq
  x, 1 > d_1 > \cdots > d_{n-1} > 0\Bigg\}
\end{align*}

Let $\mathcal{R}$ be the region with positive increasing coordinates
defined as greater in the $i$th coordinate than $\mathcal{H}_i$.  

As in the even case
\begin{lemma}\label{odd1log}
Suppose $x>0$ and 
\[
s = \{1, x, x^{1+d_1}, x^{1+d_1+d_2}, \cdots,
x^{1+d_1+\cdots+d_n}, x^{1+d_1+\cdots+d_{n-1}}, \cdots , x, 1\}.
\]
  Then
$s$ is $1$-logconcave iff $1 \geq d_1 \geq \cdots \geq d_n \geq 0$
with strict inequalities in the logconcavity condition iff $1 > d_1 > \cdots > d_n > 0$.
\end{lemma}

\begin{proof}
  The only case which differs from the proof of Lemma \ref{1log} is
  the $n$th.  $x^{2+2d_1+\cdots+2d_n} \geq x^{2+2d_1+\cdots+2d_{n-1}}
  \Leftrightarrow x^{2d_n} \geq 1 \Leftrightarrow d_n \geq 0$ and
  likewise with strict inequalities.
\end{proof}

\begin{lemma}
  $\mathcal{R}$ is nonempty and unbounded.
\end{lemma}

\begin{proof}
The proof begins as before, but with $a > (1+\sqrt{5})C^{T(n-1)-T(n)}/2$. We only need consider $\mathcal{H}_n$.
For $\mathcal{H}_n$ choose $x$ and $d_i$ as before to match
all the coordinates of $s$ other than the $n$th.  Then $x_n^2 \geq
x_{n-1}^2$ giving
\begin{align*}
  C^{T(n)}a^{n+1}x_{n} 
  & \geq
  C^{T(n)}a^{n+1}x_{n-1}
  \\
  & = aC^{T(n)-T(n-1)} x^{1+d_1+\cdots +d_{n-1}} \\
  & > \frac{1+\sqrt{5}}{2}x^{1+d_1+\cdots +d_{n-1}}
\end{align*}
So $s$ is also on the correct side of $\mathcal{H}_n$.  Consequently
$s$ is in $\mathcal{R}$.  So we see that $\mathcal{R}$ is nonempty,
and, by the freedom to increase $a$, is unbounded.
\end{proof}

\begin{lemma}\label{oddHmonotone}
  Each of the $\mathcal{H}_j$ is $j$-monotone.
\end{lemma}

\begin{proof}
  We only need to consider $\mathcal{H}_n$, in which $x_n$ is
  determined by $x_{n-1}$ and increases when 
  $x_{n-1}$ increases, so $\mathcal{H}_n$ is $n$-monotone.  
\end{proof}

\begin{theorem}
  Any sequence in $\mathcal{R}$ is $\infty$-logconcave.
\end{theorem}

\begin{proof}
  Again, in view of the even case, we only need to check
  $\mathcal{H}_n$.  Using notation from the even case choose $x$,
  $\epsilon>0$, and 
  the $d_i$, such that
  \[
    s = \Bigg\{1, x, \ldots, x^{1+d_1+\cdots+d_{n-1}},
    \frac{1+\sqrt{5}}{2}x^{1+d_1+\cdots+d_{n-1}} + \epsilon, \ldots, 1\Bigg\}
   \]
  Iterate
  to get
  \[
  \begin{split}
    \mathcal{L}(s) = \Bigg\{1, x^2-x^{1+d_1}, \ldots,
    x^{2+2d_1+\cdots+2d_{n-1}} -
    \left(\frac{1+\sqrt{5}}{2}\right)x^{2+2d_1+\cdots+d_{n-1}} - x^{1+d_1+\cdots+d_{n-2}}\epsilon, \\
    \left(\left(\frac{1+\sqrt{5}}{2}\right)^2 -
    1\right)x^{2+2d_1+\cdots+2d_{n-1}} + (1+\sqrt{5})\epsilon
    x^{1+d_1+\cdots+d_{n-1}} + \epsilon^2, \ldots, 1\Bigg\}
  \end{split}
  \]
  which by Lemmas \ref{sides} and \ref{oddHmonotone} shows
  that $\mathcal{L}(s)$ is on the same side of $\mathcal{H}_n$ as $s$
  is.

  Consequently $\mathcal{L}(s)$ is in $\mathcal{R}$.  Since
  $\mathcal{R}$ is a subregion of the region of $\mathbb{R}^{n+1}$
  with positive coordinates, this implies that any sequence in
  $\mathcal{R}$ is $\infty$-logconcave.
\end{proof}


\begin{thebibliography}{9}
\bibitem{bomo}
  Boros, G. and Moll, V.: {\em Irresistible Integrals: Symbolics,
    Analysis and Experiments in the Evaluation of
    Integrals}. Cambridge University Press, 2003.
\bibitem{br}
  Brenti, Francesco: \emph{Log-concave and Unimodal sequences in
    Algebra, Combinatorics, and Geometry: an update}.  Contemporary
  Math.,  178 (1994), 71-89. 
\bibitem{kp}
  Kauers, Manuel and Paule, Peter: \emph{A Computer Proof of Moll's
    Log-Concavity Conjecture}. SFB F13. Technical report no. 2006-15,
  Altenbergerstrasse 69, 2006.
\bibitem{m}
  Moll, Victor H.: \emph{Combinatorial sequences arising from a
    rational integral}. Online Journal of Analytic Combinatorics,
  \textbf{2} (2007), \#4.

\end{thebibliography}
\end{document}